\documentclass{amsart}
\usepackage{amssymb,latexsym,enumerate,amscd}
\newtheorem{Thm}{Theorem}[section]
\newtheorem{Prop}[Thm]{Proposition}

\newtheorem{Prob}[Thm]{Problem}

\theoremstyle{remark}

\theoremstyle{definition}

\newtheorem*{Not*}{Notation}
\newtheorem*{Def}{Definition}
\numberwithin{equation}{section}

\DeclareMathOperator{\rk}{rk}
\DeclareMathOperator{\End}{End}
\DeclareMathOperator{\R}{R}
\DeclareMathOperator{\Mat}{M}
\DeclareMathOperator{\Jac}{Jac}
\DeclareMathOperator{\diag}{diag}
\begin{document}

\date{%
Sun Jul 25 21:44:08 EDT 2004}

\title[Noncommutative localization in group rings]
{Noncommutative localization in group rings}

\author[P. A. Linnell]{Peter A. Linnell}
\address{Department of Mathematics \\
Virginia Tech \\
Blacksburg \\
VA 24061-0123 \\
USA}
\email{linnell@math.vt.edu}
\urladdr{http://www.math.vt.edu/people/linnell/}

\begin{abstract}
This paper will briefly survey some recent methods of localization
in group rings, which work in more general contexts than the
classical Ore localization.  In particular the Cohn localization
using matrices will be described, but other methods will also be
considered.
\end{abstract}



\maketitle

\section{Introduction} \label{Sintroduction}

Let $R$ be a commutative ring and let $S = \{s \in R \mid sr \ne 0$
for all $r \in R \setminus 0\}$, the set of non-zerodivisors of $R$.
Then, as in the same manner one constructs $\mathbb {Q}$
from $\mathbb {Z}$, we can form the quotient ring
$RS^{-1}$ which consists of elements of the form $r/s$ with $r \in R$
and $s \in S$, and in which $r_1/s_1 = r_2/s_2$ if and only if
$r_1s_2 = s_1r_2$.  We can consider $R$ as a
subring $RS^{-1}$ by identifying $r \in R$ with $r/1 \in RS^{-1}$.
Then $RS^{-1}$ is a ring containing $R$ with the property that every
element is either a zerodivisor or invertible.  Furthermore, every
element of $RS^{-1}$ can be written in the form $rs^{-1}$ with $r \in
R$ and $s \in S$ (though not uniquely so).  In the case $R$ is an
integral domain, then $RS^{-1}$ will be a field and will be generated
as a field by $R$ (i.e.\ if $K$ is a subfield of $RS^{-1}$ containing
$R$, then $K = RS^{-1}$).  Moreover if $K$ is another field
containing $R$ which is generated by $R$,
then $K$ is isomorphic to $RS^{-1}$ and in fact there
is a ring isomorphism $RS^{-1} \to K$ which is the identity on $R$.

The question we will be concerned with
here is what one can do with a noncommutative ring $R$;
certainly many of the above results do not hold in general.  In
particular, Malcev \cite{Malcev37} constructed domains which are not
embeddable in division rings.  We will concentrate on the case when
our ring is a crossed product $k*G$, where $k$ is a division ring
and $G$ is a
group \cite{Passman89}, and in particular when the crossed product
is the group ring $kG$ with $k$ a field.
A field will always mean a commutative
field, and we shall use the terminology ``division ring" for the
noncommutative case.  Though our main interest is in group rings,
often it is a trivial matter to extend results to crossed products.
This has the advantage of facilitating induction arguments,
because if $H \lhd G$ and $k*G$ is a crossed product, then $k*G$ can
also be viewed as a crossed product $(k*H)*(G/H)$
\cite[p.~2]{Passman89}.

\section{Ore Localization}
\label{SOre}

We shall briefly recall the definition of a crossed product, and also
establish some notational conventions for this paper.  Let $R$ be a
ring with a 1 and let $G$ be a group.  Then a crossed product of $G$
over $R$ is an associative ring $R*G$ which is also a free left
$R$-module with basis $\{\bar{g} \mid g \in G\}$.  Multiplication is
given by $\bar{x}\bar{y} = \tau(x,y)\overline{xy}$ where
$\tau(x,y)$ is a unit of $R$ for all $x,y \in G$.  Furthermore we
assume that $\bar{1}$ is the identity of $R*G$, and we identify
$R$ with $R\bar{1}$ via $r \mapsto r\bar{1}$.  Finally $\bar{x} r =
r^{\sigma(x)}\bar{x}$ where $\sigma(x)$ is an automorphism of $R$ for
all $x \in G$; see \cite[p.~2]{Passman89} for further details.

We shall assume that all rings have a 1, subrings have the same
1, and ring homomorphisms preserve the 1.
We say that the element $s$ of $R$ is a non-zerodivisor (sometimes
called a regular element) if $sr \ne 0 \ne rs$ whenever $0 \ne r \in
R$; otherwise $s$ is called a zerodivisor.
Let $S$ denote the set of non-zerodivisors of the ring $R$.
The simplest extension to noncommutative rings is when the ring
$R$ satisfies the \emph{right Ore condition}, that is given
$r \in R$ and $s \in S$, then there
exists $r_1 \in R$ and $s_1 \in S$ such that $rs_1 = sr_1$.  In this
situation one can form the Ore localization $RS^{-1}$, which in the
same way as above consists of elements of the form $\{rs^{-1}
\mid r \in R,\,s \in S\}$.  If $s_1s = s_2r$, then $r_1s_1^{-1} =
r_2s_2^{-1}$ if and only if $r_1s = r_2r$; this does not depend on
the choice of $r$ and $s$.  To define addition in
$RS^{-1}$, note that any two elements can be written in the form
$r_1s^{-1}, r_2s^{-1}$ (i.e.\ have the same common denominator), and
then we set $r_1s^{-1} + r_2s^{-1} = (r_1 + r_2)s^{-1}$.  To define
multiplication, if $s_1r = r_2s$, we set
$(r_1s_1^{-1})(r_2s_2^{-1}) = r_1r(s_2s)^{-1}$.  Then $RS^{-1}$ is a
ring with $1 = 11^{-1}$ and $0 = 0 1^{-1}$, and $\{r1^{-1} \mid r \in
R\}$ is a subring isomorphic to $R$ via the map $r \mapsto
r1^{-1}$.  Furthermore $RS^{-1}$ has the following properties:
\begin{itemize}
\item
Every element of $S$ is invertible in $RS^{-1}$.

\item
Every element of $RS^{-1}$ is either invertible or a zerodivisor.

\item
If $\theta \colon R \to K$ is a ring homomorphism such that
$\theta s$ is invertible for all $s \in S$, then there is a unique
ring homomorphism $\theta' \colon RS^{-1} \to K$ such that
$\theta' (r1^{-1}) = \theta r$ for all $r \in R$; in other words,
$\theta$ can be extended in a unique way to $RS^{-1}$.

\item
$RS^{-1}$ is a flat left $R$-module
\cite[Proposition II.3.5]{Stenstrom75}.
\end{itemize}
Of course one also has the left Ore condition, which means that given
$r \in R$ and $s\in S$, one can find $r_1 \in R$ and $s_1 \in S$
such that $s_1r = r_1s$, and then one can form the ring $S^{-1}R$,
which consists of elements of the form $s^{-1}r$ with $s \in S$ and
$r \in R$.  However in the case of the group ring $kG$ for a field
$k$ and group $G$, they are equivalent by using the involution
on $kG$ induced by $g \mapsto g^{-1}$ for $g \in G$.  When a ring
satisfies both the left and right Ore condition, then the rings
$S^{-1}R$ and $RS^{-1}$ are isomorphic, and can be identified.
In this situation, we say that $RS^{-1}$ is
a classical ring of quotients for $R$.  When $R$ is a domain, a
classical ring of quotients will be a division ring.  On the other
hand if already every element of $R$ is either invertible or a
zerodivisor, then $R$ is its own classical quotient ring.
For more information on Ore localization, see \cite[\S
9]{GoodearlWarfield89}.

\begin{Prob} \label{Pprob1}
Let $k$ be a field.  For which groups $G$ does $kG$ have a classical
quotient ring?
\end{Prob}

One could ask more generally given a division ring $D$, for which
groups $G$ does a crossed product $D*G$ always
have a classical quotient ring?  We have put in the ``always"
because $D$ and $G$ do not determine a crossed product
$D*G$.  One could equally consider the same question with ``always"
replaced by ``never".

For a nonnegative integer $n$, let $F_n$ denote the free group on $n$
generators, which is nonabelian for $n \ge 2$.
If $G$ is abelian in Problem \ref{Pprob1}, then $kG$ certainly has a
classical quotient ring because $kG$ is commutative in this case.  On
the other hand if $G$ has a subgroup isomorphic to $F_2$, then $D*G$
cannot have a classical quotient ring.  We give an elementary proof of
this well-known statement, which is based on \cite[Theorem
1]{Lewin70}.

\begin{Prop} \label{PLewin}
Let $G$ be a group which has a subgroup isomorphic to the free group
$F_2$ on two generators, let $D$ be a division ring, and let $D*G$ be
a crossed product.  Then $D*G$ does not satisfy the right Ore
condition, and in particular does not have a classical quotient ring.
\end{Prop}
\begin{proof}
First suppose $G$ is free on $a,b$.  We prove that
$(\bar{a}-1)D*G \cap (\bar{b}-1)D*G
= 0$.  Write $A = \langle a \rangle$ and $B = \langle b \rangle$.
Suppose $\alpha \in (\bar{a}-1)D*G \cap (\bar{b}-1)D*G$.
Then we may write
\begin{equation} \label{Eore}
\alpha =
\sum_i (u_i - 1)\overline{x_i}d_i
= \sum_i(v_i - 1)\overline{y_i}e_i
\end{equation}
where $u_i = \bar{a}^{q(i)}$ for some ${q(i)} \in \mathbb {Z}$, $v_i =
\bar{b}^{r(i)}$ for some ${r(i)} \in \mathbb {Z}$,
$d_i,e_i \in D$ and $x_i,y_i \in G$.  The general element
$g$ of $G$ can be written in a unique way $g_1\dots g_l$, where
the $g_i$ are alternately in $A$ and $B$,
and $g_i \ne 1$ for all $i$; we shall define the length
$\lambda (g)$ of $g$ to be $l$.  Of course $\lambda (1) = 0$.  Let
$L$ be the maximum of all $\lambda (x_i),\lambda (y_i)$,
let $s$ denote the number of
$x_i$ with $\lambda(x_i) = L$, and let $t$ denote the
number of $y_i$ with $\lambda(y_i) = L$.
We shall use induction on $L$ and then on $s+t$, to show
that $\alpha = 0$.  If
$L = 0$, then $x_i,y_i = 1$ for all $i$ and the result is obvious.
If $L > 0$, then without loss of generality, we may assume that
$s>0$.  Suppose $\lambda(x_i) = L$ and
$x_i$ starts with an element from $A$,
so $x_i = a^ph$ where $0 \ne p \in \mathbb {Z}$
and $\lambda(h) = L-1$.  Then
\[
(u_i - 1)\overline{x_i}d_i =
(\bar{a}^{q(i)} - 1)\bar{a}^p\bar{h}dd_i
= (\bar{a}^{q(i)+p}-1)\bar{h}dd_i
- (\bar{a}^p - 1) \bar{h}dd_i
\]
for some $d \in D$.  This means that we have found an expression for
$\alpha$ with smaller $s+t$, so all the $x_i$ with $\lambda(x_i) = L$
start with an element from $B$.  Therefore if $\beta = \sum_i u_i
\overline{x_i}d_i$ where the sum is over all $i$ such that
$\lambda(x_i) = L$, then each $x_i$ starts with an element of $B$ and
hence $\lambda (a^{q(i)}x_i) = L+1$.
We now see from \eqref{Eore} that $\beta = 0$.  Since $s>0$ by
assumption, the expression for $\beta$ above is nontrivial and
therefore there exists $i \ne j$ such that $a^{q(i)}x_i =
a^{q(j)}x_j$.  This forces $q(i) = q(j)$ and $x_i = x_j$.
Thus $u_i = u_j$ and we may replace
$(u_i - 1)\overline{x_i}d_i + (u_j - 1)\overline{x_j}d_j$ with
$(u_i - 1)\overline{x_i}(d_i+d_j)$, thereby reducing $s$ by 1 and the
proof that $(\bar{a}-1)D*G \cap (\bar{b}-1)D*G = 0$ is complete.

In general, suppose $G$ has a subgroup $H$ which is free on the
elements $x,y$.  Then the above shows that $(x-1)D*H \cap (y-1)D*H =
0$, and it follows that $(x-1)D*G \cap (y-1)D*G = 0$.
Since $x-1$ and $y-1$ are non-zerodivisors in $D*G$, it follows that
$D*G$ does not have the right Ore property.
\end{proof}

Recall that the class of elementary amenable groups is the smallest
class of groups which contains all finite groups and the infinite
cyclic group $\mathbb {Z}$, and is closed under taking group
extensions and directed unions.  It is not difficult to show that the
class of elementary amenable groups is closed under taking subgroups
and quotient groups,
and contains all solvable-by-finite groups.  Moreover every
elementary amenable group is amenable, but $F_2$ is not amenable.
Thus any group which has a subgroup
isomorphic to $F_2$ is not elementary amenable.  Also
Thompson's group $F$ \cite[Theorem 4.10]{CFP96}
and the Gupta-Sidki group \cite{Gupta89}
are not elementary amenable even though they do not contain $F_2$.
The Gupta-Sidki has sub-exponential growth \cite{FabrykowskiGupta85}
and is therefore amenable \cite[Proposition 6.8]{Paterson88}.
The following result follows from \cite[Theorem 1.2]{KLM88}

\begin{Thm} \label{TKLM}
Let $G$ be an elementary amenable group, let $D$ be a division
ring, and let $D*G$ be a crossed product.
If the finite subgroups of $G$ have bounded order, then $D*G$ has a
classical ring of quotients.
\end{Thm}

It would seem plausible that Theorem
\ref{TKLM} would remain true without the
hypothesis that the finite subgroups have bounded order.  After all,
if $G$ is a locally finite group and $k$ is a field, then $kG$ is a
classical quotient ring for itself.  However the lamplighter group,
which we now describe, yields a counterexample.  If $A,C$ are groups,
then $A \wr C$ will indicate the Wreath product with base group $B :=
A^{|C|}$, the direct sum of $|C|$ copies of $A$.  Thus $B$ is a
normal subgroup of $A \wr C$ with corresponding quotient group
isomorphic to $C$, and $C$ permutes the $|C|$ copies of $A$
regularly.  The case $A = \mathbb {Z}/2\mathbb {Z}$
and $C = \mathbb {Z}$ is often called the lamplighter group.
Then \cite[Theorem 2]{LLS} is
\begin{Thm} \label{TLLS}
Let $H\ne 1$ be a finite group, let $k$ be a field, and let $G$
be a group containing $H \wr \mathbb {Z}$.
Then $kG$ does not have a classical ring of quotients.
\end{Thm}

Thus we have the following problem.
\begin{Prob}
Let $k$ be a field.  Classify the elementary amenable groups $G$
for which $kG$ has a classical ring of quotients.  If $H \leqslant G$
and $kG$ has a classical ring of quotients, does $kH$ also have a
classical ring of quotients?
\end{Prob}

The obstacle to preventing a classical quotient ring in the case of
elementary amenable groups is the finite subgroups having unbounded
order, so let us consider the case of torsion-free groups.  In this
situation it is unknown whether $kG$ is a domain, so let us assume
that this is the case.  Then we have the following result of
Tamari \cite{Tamari54}; see \cite[Theorem 6.3]{DLMSY03}, also
\cite[Example 8.16]{Lueck02}, for a proof.

\begin{Thm} \label{TTamari}
Let $G$ be an amenable group, let $D$ be a division ring,
and let $D*G$ be a crossed product which is a domain.  Then $D*G$ has
a classical ring of quotients which is a division ring.
\end{Thm}

What about torsion-free groups which do not contain $F_2$, yet are
not amenable?  Given such a group $G$ and a division ring
$D$, it is unknown
whether a crossed product $D*G$ has a classical quotient ring.
Thompson's group $F$ is orderable \cite[Theorem 4.11]{CFP96}; this
means that it has a total order $\le$ which is left and right
invariant, so if $a \le b$ and $g \in F$, then $ga \le gb$ and $ag
\le bg$.  Therefore if $D$ is a division ring and $D*F$ is a crossed
product, then by the Malcev-Neumann construction \cite[Corollary
8.7.6]{Cohn85} the power series ring $D((F))$ consisting of elements
with well-ordered support is a division ring.
It is still unknown whether
Thompson's group is amenable.  We state the following problem.
\begin{Prob} \label{PThompson}
Let $F$ denote Thompson's orderable group and let $D$ be a division
ring.  Does $D*F$ have a classical ring of quotients?
\end{Prob}
If the answer is negative, then Theorem \ref{TTamari}
would tell us that Thompson's group is not amenable.  Since
Thompson's group seems to be right on the borderline between
amenability and nonamenability, one would expect the answer to be in
the affirmative.

\section{Cohn's Theory} \label{SCohn}

What happens when the ring $R$ does not have the Ore condition, in
other words $R$ does not have a classical ring of quotients?
Trying to form a ring from $R$ by inverting the non-zerodivisors of
$R$ does not seem very useful.
The key idea here is due to Paul Cohn; instead of trying to invert
just elements, one inverts matrices instead.  Suppose $\Sigma$ is
any set of matrices over $R$ (not necessarily square, though in
practice $\Sigma$ will consist only of square matrices)
and $\theta \colon R \to S$ is a ring
homomorphism.  If $M$ is a matrix with entries $m_{ij} \in R$, then
$\theta M$ will indicate the matrix over $S$ which has entries
$\theta(m_{ij})$.  We say that $\theta$ is $\Sigma$-inverting if
$\theta M$ is invertible over $S$ for all $M \in \Sigma$.
We can now define the universal localization of $R$ with respect to
$\Sigma$, which consists of a ring $R$ and a \emph{universal}
$\Sigma$-inverting ring homomorphism
$\lambda \colon R \to R_{\Sigma}$.  This means that given any other
$\Sigma$-inverting homomorphism $\theta \colon R \to S$, then there
is a unique ring homomorphism $\phi \colon R_{\Sigma} \to S$ such
that $\theta = \phi\lambda$.  The ring $R_{\Sigma}$ always exists by
\cite[Theorem 7.2.1]{Cohn85}, and by the universal property is unique
up to isomorphism.  Furthermore $\lambda$ is injective if and only if
$R$ can be embedded in a ring over which all the matrices in $\Sigma$
become invertible.

A related concept is the $\Sigma$-rational closure.  Given a set of
matrices $\Sigma$ over $R$ and a $\Sigma$-inverting ring
homomorphism $\theta \colon R\to S$,
the $\Sigma$-rational closure $R_{\Sigma}(S)$ of $R$ in $S$ consists
of all entries of inverses of matrices in $\theta(\Sigma)$.
In general $R_{\Sigma}(S)$ will not be a subring of $S$.  We say
that $\Sigma$ is upper multiplicative if given $A,B \in \Sigma$,
then
$\begin{pmatrix}
A&C\\
0&B
\end{pmatrix} \in \Sigma$
for any matrix $C$ of the appropriate size.  If in addition permuting
the rows and columns of a matrix in $\Sigma$ leaves it in $\Sigma$,
then we say that $\Sigma$ is multiplicative.

Suppose now that $\Sigma$ is a set of matrices over $R$ and
$\theta \colon R\to S$ is a $\Sigma$-inverting ring homomorphism.
If $\Sigma$ is upper multiplicative,
then $R_{\Sigma}(S)$ is a subring of $S$ \cite[Theorem
7.1.2]{Cohn85}.  Also if $\Phi$ is the set of matrices over $R$ whose
image under $\theta$ becomes invertible over $S$, then $\Phi$ is
multiplicative \cite[Proposition 7.1.1]{Cohn85}.  In this situation
we call $R_{\Phi}(S)$ the rational closure $\mathcal{R}_S(R)$
of $R$ in $S$.  By the universal property of $R_{\Phi}$, there is a
ring homomorphism $R_{\Phi} \to R_{\Phi}(S) = \mathcal{R}_S(R)$
which is
surjective.  A very useful tool is the following consequence of
\cite[Proposition 7.1.3]{Cohn85}, which we shall call Cramer's rule;
we shall let $\Mat_n(R)$ denote the $n \times n$ matrices over $R$.
\begin{Prop} \label{PCramer}
Let $\Sigma$ be an upper multiplicative set of matrices of $R$ and
let $\theta \colon R\to S$ be a $\Sigma$-inverting
ring homomorphism.  If
$p \in R_{\Sigma}(S)$, then $p$ is stably associated to a matrix with
entries in $\theta(R)$.  This means that there exists a positive
integer $n$ and invertible matrices $A,B \in \Mat_n(S)$ such that
$A\diag(p,1,\dots,1)B \in \Mat_n(\theta R)$.
\end{Prop}

Given a ring homomorphism $\theta \colon R \to S$ and
an upper multiplicative set of matrices $\Sigma$ of $R$,
the natural epimorphism $R_{\Sigma} \twoheadrightarrow
R_{\Sigma}(S)$ will in general not be isomorphism, even if $\theta$
is injective, but there are interesting situations where it is; we
describe one of them.  Let $k$
be a PID (principal ideal domain), let $X$ be a set, let
$k\langle X\rangle$ denote the free algebra on $X$, let
$k\langle\langle X\rangle\rangle$ denote the noncommutative power
series ring on $X$, and let $\Lambda$ denote the subring of
$k\langle\langle X\rangle\rangle$ generated by $k\langle X\rangle$
and $\{(1+x)^{-1} \mid x \in X\}$.  Then $\Lambda \cong kF$ where $F$
denotes the free group on $X$ \cite[p.~529]{Cohn85}.  Let $\Sigma$
consist of those square matrices over $\Lambda$ with constant term
invertible over $k$, and let $\Sigma' = \Sigma\cap k\langle X\rangle$.
If we identify $\Lambda$ with $kF$ by the above isomorphism, then
$\Sigma$ consist of those matrices over $kF$ which become invertible
under the augmentation map $kF \to k$.  Since $\Sigma$ and $\Sigma'$
are precisely the matrices over $\Lambda$ and $k\langle X\rangle$
which become invertible over $k\langle\langle X\rangle\rangle$
respectively, we see that
$k\langle X\rangle_{\Sigma'} (k\langle\langle X\rangle\rangle) =
\Lambda_{\Sigma} (k\langle\langle X\rangle\rangle) =
\mathcal{R}_{k\langle\langle X\rangle\rangle}(k\langle X\rangle) =
\mathcal{R}_{k\langle\langle X\rangle\rangle}(\Lambda)$.
By universal properties, we have a sequence of natural maps
\[
\Lambda \overset{\alpha}{\longrightarrow}
k\langle X\rangle_{\Sigma'} \overset{\beta}{\longrightarrow}
\Lambda_{\Sigma} \overset{\gamma}{\longrightarrow}
k\langle X\rangle_{\Sigma'} (k\langle\langle X\rangle\rangle).
\]
The map $\gamma\beta$ is an isomorphism by
\cite[Theorem 24]{DicksSontag78}.  Therefore the image under
$\alpha$ of every matrix in $\Sigma$ becomes invertible in $k\langle
X\rangle_{\Sigma'}$, hence there is a natural map $\phi \colon
\Lambda_{\Sigma} \to k\langle X\rangle_{\Sigma'}$ such that
$\beta\phi$ and $\phi\beta$ are the identity maps.  We deduce
that $\gamma$ is also an isomorphism.  It would be interesting to know
if $\gamma$ remains an isomorphism if $k$ is assumed to be only an
integral domain.  We state the following problem.

\begin{Prob}
Let $X$ be a set, let $F$ denote the free group on $X$, and
let $k$ be an integral domain.  Define a $k$-algebra monomorphism
$\theta \colon kF \to k\langle\langle X\rangle\rangle$ by $\theta (a)
= a$ for $a \in k$ and $\theta (x) = 1+x$ for $x \in X$, let $\Sigma$
be the set of matrices over $kF$ which become invertible over
$k\langle\langle X\rangle\rangle$ via $\theta$, and let $\phi \colon
kF_{\Sigma} \to k\langle\langle X\rangle\rangle$ be the uniquely
defined associated ring homomorphism.  Determine when $\phi$ is
injective.
\end{Prob}

If $R$ is a subring
of the ring $T$, then we define the division closure
$\mathcal{D}_T(R)$ of $R$ in $T$ to
be the smallest subring $\mathcal{D}_T(R)$ of $T$ containing $R$
which is closed under taking inverses,
i.e.\ $x \in \mathcal{D}_T(R)$ and $x^{-1} \in T$ implies
$x^{-1} \in \mathcal{D}_T(R)$.
In general $\mathcal{D}_T(R) \subseteq \mathcal{R}_T(R)$,
i.e.\ the division closure is contained in
the rational closure \cite[Exercise 7.1.1]{Cohn85}.  However if $T$
is a division ring, then the rational closure is a division ring and
is equal to the division closure.

It is clear that taking the division closure is an idempotent
operation; in other words $\mathcal{D}_T(\mathcal{D}_T(R)) =
\mathcal{D}_T(R)$.  It is also true that taking the rational closure
is an idempotent operation; we sketch the proof below.

\begin{Prop}
Let $R$ be a subring of the ring $T$ and assume that $R$ and $T$ have
the same 1.  Then $\mathcal{R}_T(\mathcal{R}_T(R)) = \mathcal
{R}_T(R)$.
\end{Prop}
\begin{proof}
Write $R' =\mathcal{R}_T(R)$ and let $M$ be a matrix over $R'$ which
is invertible over $T$; we need to prove that
all the entries of $M^{-1}$ are in $R'$.  We may assume that $M \in
\Mat_d(R')$ for some positive integer $d$.  Cramer's
rule, Proposition \ref{PCramer}, applied to the inclusion
$\Mat_d(R) \to \Mat_d(R')$ tells us that $M$ is
stably associated to a matrix with entries in $\Mat_d(R)$.  This
means that for some positive integer $e$, there exists a matrix
$L\in\Mat_e(\Mat_d(R)) = \Mat_{de}(R)$ of
the form $\diag(M,1,\dots,1)$ and invertible matrices $A,B \in
\Mat_{de}(R')$ such that $ALB$ is a matrix $X \in \Mat_{de}(R)$.

Since $A,L,B$ are all invertible in
$\Mat_{de}(T)$, we see that
$X^{-1}$ has (by definition of rational closure) all its entries in
$\Mat_{de}(R')$.  But $L^{-1} = BX^{-1}A$, which shows that $L^{-1}
\in \Mat_{de}(R')$.  Therefore $M^{-1} \in \Mat_d(R')$ as required.
\end{proof}
We also have the following useful result.

\begin{Prop} \label{Pmatrix}
Let $n$ be a positive integer, let $R$ be a subring of the ring
$T$, and assume that $R$ and $T$ have the same 1.  Then
$\mathcal{R}_{\Mat_n(T)}(\Mat_n(R)) = \Mat_n(\mathcal{R}_T(R))$.
\end{Prop}
\begin{proof}
Write $R' = \mathcal{R}_T(R)$ and $S = \Mat_n(T)$.
Suppose $M \in \mathcal{R}_{S}(\Mat_n(R))$.  Then $M$
appears as an entry of $A^{-1}$, where $A \in \Mat_d(\Mat_n(R))$
for some positive integer $d$ is
invertible in $\Mat_d(S)$.  By definition all the entries of $A^{-1}$
(when viewed as a matrix in $\Mat_{dn}(T)$)
are in $R'$, which shows that $M \in \Mat_n(R')$.

Now let $M \in \Mat_n(R')$.  We want to show that $M\in
\mathcal{R}_S(\Mat_n(R))$.  Since
$\mathcal{R}_S(\Mat_n(R))$ is a ring, it is
closed under addition, so
we may assume that $M$ has exactly one nonzero entry.  Let $a$ be
this entry.  Then $a$ appears
as an entry of $A^{-1}$ where $A$ is an invertible matrix in
$\Mat_m(R)$ for some positive
integer $m$ which is a multiple of $n$.
By permuting the rows and columns, we may assume that
$a$ is the $(1,1)$-entry.  Now
form the $p \times p$ matrix $B = \diag(1,\dots,1,A,1,\dots,1)$, so
that the $(1,1)$-entry of $A$ is in the $(n,n)$-entry of
$B$ (thus there are $n-1$ ones on the main diagonal
and then $A$) and $m$ divides $p$.  By considering $B^{-1}$, we
see that $\diag(1,\dots,1,a) \in \mathcal{R}_S(\Mat_n(R))$.
Since $\diag(1,\dots,1,0) \in \mathcal{R}_S(\Mat_n(R))$,
it follows that $\diag(0,\dots,0,a) \in
\mathcal{R}_S(\Mat_n(R))$.  By permuting the
rows and columns, we conclude that
$M \in \mathcal{R}_S(\Mat_n(R))$.
\end{proof}

When one performs a localization, it would be good to end up with a
local ring.  We now describe a result of
Sheiham \cite[\S 2]{Sheiham02} which
shows that this is often the case.
For any ring $R$, we let $\Jac(R)$ indicate the Jacobson
radical of $R$.  Let $\theta \colon R \to S$ be a ring
homomorphism, let
$\Sigma$ denote the set of all matrices $A$ over $R$ with the
property that $\theta (A)$ is an invertible matrix over $S$, and let
$\lambda \colon R \to R_{\Sigma}$ denote the associated map.  Then we
have a ring homomorphism $\phi \colon R_{\Sigma} \to S$ such that
$\theta = \phi \lambda$, and Sheiham's result is
\begin{Thm}
If $S$ is a local ring, then
$\phi^{-1}\Jac(S) = \Jac(R_{\Sigma})$
\end{Thm}
Thus in particular if $S$ is a division ring, then $R_{\Sigma}$ is a
local ring.

\section{Uniqueness of Division Closure
and Unbounded Operators} \label{Suniqueness}

If $R$ is a domain and $D$ is a division ring containing $R$ such
that $\mathcal{D}_D(R) = D$
(i.e.\ $R$ generates $D$ as a field), then
we say that $D$ is a division ring of fractions for $R$.  If $R$ is
an integral domain and $D,E$ are division rings of fractions for $R$,
then $D$ and $E$ are fields and are just the Ore localizations of $R$
with respect to the nonzero elements of $R$.  In this case there
exists a unique isomorphism $D \to E$ which is the identity on $R$.
Furthermore any automorphism of $D$ can be extended
to an automorphism of $R$.

When $D$ and $E$ are not commutative, i.e.\ it is
only assumed that they are division rings, then this is not the case;
in fact $D$ and $E$ may not be isomorphic even just as rings.
Therefore we would like to have a criterion for when two such
division rings are isomorphic, and also a criterion for the closely
related property of when an automorphism of $R$ can be extended to an
automorphism of $D$.

Consider now the complex group algebra $R=\mathbb {C}G$.  Here we may
embed $\mathbb {C}G$ into the ring of unbounded operators
$\mathcal{U}(G)$ on $L^2(G)$
affiliated to $\mathbb {C}G$; see e.g.\ \cite[\S 8]{Linnell98} or
\cite[\S 8]{Lueck02}.  We briefly recall the construction and state
some of the properties.  Let $L^2(G)$ denote the Hilbert space with
Hilbert basis the elements of $G$; thus $L^2(G)$ consists of all
square summable formal sums $\sum_{g \in G} a_g g$ with
$a_g \in \mathbb {C}$ and inner product $\langle \sum_g
a_gg,\sum_h b_hh\rangle = \sum_{g,h} a_g \overline{b_h}$.
We have a left and right action of $G$ on $L^2(G)$ defined by the
formulae $\sum_h a_hh \mapsto \sum_h a_hgh$ and
$\sum_h a_hh \mapsto \sum_h a_hhg$ for $g \in G$.  It follows that
$\mathbb {C}G$ acts faithfully
as bounded linear operators on the left of
$L^2(G)$, in other words we may consider $\mathbb {C}G$ as a subspace
of $\mathcal{B}(L^2(G))$, the bounded linear operators on $L^2(G)$.
The weak closure of $\mathbb {C}G$ in $\mathcal{B}(L^2(G))$ is the
group von Neumann algebra $\mathcal{N}(G)$ of $G$, and the unbounded
operators affiliated to $G$, denoted $\mathcal{U}(G)$, are those
closed densely defined unbounded operators which commute
with the right action of $G$.  We have a natural injective
$\mathbb{C}$-linear map $\mathcal{N}(G) \to L^2(G)$ defined by
$\theta \mapsto \theta 1$ (where 1 denotes the element $1_1$ of
$L^2(G))$, so we may identify $\mathcal{N}(G)$ with a subspace of
$L^2(G)$.  When $H \leqslant G$, we may consider $L^2(H)$ as a
subspace of $L^2(G)$ and using the above identification, we may
consider $\mathcal{N}(H)$ as a subring of $\mathcal{N}(G)$.
Also given $\alpha \in L^2(G)$, we can define a $\mathbb {C}$-linear
map $\hat{\alpha} \colon \mathbb {C}G \to L^2(G)$ by
$\hat{\alpha}(\beta) = \alpha\beta$ for $\beta \in L^2(G)$.  Since
$\mathbb {C}G$ is a dense linear subspace of $L^2(G)$, it yields a
densely defined unbounded operator on $L^2(G)$ which commutes with
the right action of $G$, and it is not difficult to see that this
defines a unique element of $\mathcal{U}(G)$, which we shall also
call $\hat{\alpha}$.  We now have $\mathcal{N}(G) \subseteq L^2(G)
\subseteq \mathcal{U}(G)$.  Obviously if $G$ is finite, then
$\mathcal{N}(G) = L^2(G) = \mathcal{U}(G)$, because all terms are
equal to $\mathbb {C}G$.  This raises the following question.
\begin{Prob}
Let $G$ be an infinite group.  Is it always the true that
$\mathcal{N}(G) \ne L^2(G) \ne \mathcal{U}(G)$?
\end{Prob}
Presumably the answer is yes, but I am not aware of any reference.
Some related information on this, as well as results on various
homological dimensions of $\mathcal{N}(G)$ and $\mathcal{U(G)}$, can
be found in \cite{Vas02}.

At this stage it is less important to understand the construction of
$\mathcal{U}(G)$ than to know its properties.  Recall that $R$ is a
von Neumann regular ring means that given $r \in R$, there exists $x
\in R$ such that $rxr=r$.  All matrix rings over a von Neumann
regular ring are also von Neumann regular \cite[Lemma
1.6]{Goodearl91},
and every element of a von
Neumann regular ring is either invertible or a zerodivisor.
We now have that
$\mathcal{U}(G)$ is a von Neumann regular ring containing
$\mathcal{N}(G)$, and is a classical ring of quotients for
$\mathcal{N}(G)$ \cite[proof of Theorem 10]{Berberian82} or
\cite[Theorem 8.22(1)]{Lueck02}.
Thus the embedding of $\mathcal{N}(H)$ in $\mathcal{N}(G)$ for $H
\leqslant G$ as described above extends to a natural embedding of
$\mathcal{U}(H)$ in $\mathcal{U}(G)$.
Also $\mathcal{U}(G)$ is rationally closed in any overing.
Furthermore $\mathcal{U}(G)$ is a self injective unit-regular
ring which is the maximal ring of quotients of $\mathcal{N}(G)$
\cite[Lemma 1, Theorems 2 and 3]{Berberian82}.
Thus we have embedded $\mathbb {C}G$ in
a ring, namely $\mathcal{U}(G)$, in which
every element is either invertible or a zerodivisor.  In fact every
element of any matrix ring over $\mathcal{U}(G)$ is either invertible
or a zerodivisor.  Of course the
same is true for any subfield $k$ of $\mathbb {C}$, that is $kG$ can
be embedded in a ring in which every element is either invertible or
a zerodivisor.  Let us write $\mathcal{D}(kG) =
\mathcal{D}_{\mathcal{U}(G)}(kG)$ and $\mathcal{R}(kG) =
\mathcal{R}_{\mathcal{U}(G)}(kG)$.  Then if $H\leqslant G$, we may
by the above identify $\mathcal{D}(kH)$ with
$\mathcal{D}_{\mathcal{U}(G)}(kH)$ and $\mathcal{R}(kH)$ with
$\mathcal{R}_{\mathcal{U}(G)}(kH)$.
More generally, we shall write $\mathcal{D}_n(kG) =
\mathcal{D}_{\Mat_n(\mathcal{U}(G))}(\Mat_n(kG))$
and $\mathcal{R}_n(kG) =
\mathcal{R}_{\Mat_n(\mathcal{U}(G))}(\Mat_n(kG))$.  Thus
$\mathcal{D}_1(kG) = \mathcal{D}(kG)$ and
$\mathcal{R}_1(kG) = \mathcal{R}(kG)$.  Also, we may
identify $\mathcal{D}_n(kH)$ with
$\mathcal{D}_{\Mat_n(\mathcal{U}(G))}(\Mat_n(kH))$
and $\mathcal{R}(kG)$ with
$\mathcal{R}_{\Mat_n(\mathcal{U}(G))}(\Mat_n(kH))$.

Often $\mathcal{R}(kG)$ is a very nice ring.  For example
when $G$ has a normal free subgroup with elementary amenable
quotient, and also the finite subgroups of $G$ have bounded order,
it follows from \cite[Theorem 1.5(ii)]{Linnell93} that
$\mathcal{R}(\mathbb {C}G)$ is a semisimple Artinian ring, i.e.\ a
finite direct sum of matrix rings over division rings.
Thus in particular every element of $\mathcal{R}(\mathbb {C}G)$
is either invertible or a zerodivisor.  We state the following
problem.

\begin{Prob} \label{Pinvertible}
Let $G$ be a group and let $k$ be a subfield of $\mathbb {C}$.  Is
every element of $\mathcal{R}_n(kG)$ either invertible or a
zerodivisor for all positive integers $n$?
Furthermore is $\mathcal{D}_n(kG) = \mathcal{R}_n(kG)$?
\end{Prob}
The answer is certainly in the affirmative if $G$ is amenable.

\begin{Prop} \label{Pelek}
Let $G$ be an amenable group, let $n$ be a positive integer, and let
$k$ be a subfield of $\mathbb {C}$.  Then every element of
$\mathcal{R}_n(kG)$ is either a zerodivisor or invertible.
Furthermore $\mathcal{D}_n(kG) = \mathcal{R}_n(kG)$.
\end{Prop}
\begin{proof}
Write $R = \mathcal{R}_n(kG)$ and let $A \in R$.
By Cramer's rule Proposition \ref{PCramer},
there is a positive integer $d$ and invertible
matrices $X,Y \in \Mat_d(R)$ such that $B:=X\diag(A,1,\dots,1)Y \in
\Mat_{dn}(kG)$.  Suppose $ZA\ne 0 \ne AZ$ whenever $0 \ne Z \in
\Mat_n(kG)$.  Then $B$ is a
non-zerodivisor in $\Mat_{dn}(kG)$.  We claim that $B$ is also a
non-zerodivisor in $\Mat_{dn}(\mathbb {C}G)$.  If our claim is false,
then either $BC=0$ or $CB=0$ for some nonzero
$C \in \Mat_{dn}(\mathbb {C}G)$.  Without loss of generality, we may
assume that $BC=0$.  Then
for some positive integer $m$, we may choose $e_1,\dots,e_m \in
\mathbb{C}$ which are linearly independent over $k$ such that we may
write $C = C_1e_1 + \dots + C_me_m$,
where $0 \ne C_i \in \Mat_{dn}(kG)$ for all $i$.
The equation $BC= 0$ now yields $BC_1 = 0$, contradicting the
fact that $B$ is a non-zerodivisor in $\Mat_{dn}(kG)$, and the claim
is established.

Now $B$ induces by left multiplication a right $\mathbb
{C}G$-monomorphism $\mathbb {C}G^{dn} \to \mathbb {C}G^{dn}$.
This in turn induces a right $\mathcal{N}(G)$-map
$\mathcal{N}(G)^{dn} \to \mathcal{N}(G)^{dn}$, and the kernel of
this map has dimension 0 by \cite[Theorem 5.1]{Lueck98}.
It now follows from the theory of \cite[\S
2]{Lueck98} that this kernel is 0, consequently $B$ is a
non-zerodivisor in $\Mat_{dn}(\mathcal{N}(G))$.  Since
$\mathcal{U}(G)$ is a classical ring of quotients for
$\mathcal{N}(G)$, we see that $B$ is invertible
in $\Mat_{dn}(\mathcal{U}(G))$ and hence $B$ is invertible in
$\Mat_d(R)$.  Therefore $A$ is invertible in $R$ and the result
follows.
\end{proof}
One could ask the following stronger problem.
\begin{Prob} \label{PvonNeumann}
Let $G$ be a group and let $k$ be a subfield of $\mathbb {C}$.  Is
$\mathcal{R}(kG)$ a von Neumann regular ring?
\end{Prob}
Since being von Neumann regular is preserved under Morita equivalence
\cite[Lemma 1.6]{Goodearl91} and $\mathcal{R}_n(kG)$ can be
identified with $\Mat_n(\mathcal{R}(kG))$ by Proposition
\ref{Pmatrix}, we see that this is equivalent to asking
whether $\mathcal{R}_n(kG)$ is a von Neumann regular ring.
Especially interesting is the case of the lamplighter group,
specifically
\begin{Prob} \label{Plamplighter}
Let $G$ denote the lamplighter group.  Is $\mathcal{R}(\mathbb {C}G)$
a von Neumann regular ring?
\end{Prob}

Suppose $H \leq G$ and $T$ is a right transversal for $H$ in $G$.
Then $\bigoplus_{t\in T} L^2(H)t$
is a dense linear subspace of $L^2(G)$, and
$\mathcal{U}(H)$ is naturally a subring of $\mathcal{U}(G)$ as
follows.  If $u \in
\mathcal{U}(H)$ is defined on the dense linear subspace $D$ of
$L^2(H)$, then we can extend $u$ to the dense linear subspace
$\bigoplus_{t\in T} Dt$ of $L^2(G)$ by the rule $u (dt) = (ud)t$ for
$t \in T$, and the resulting unbounded
operator commutes with the right action of $G$.  It is not difficult
to show that $u \in \mathcal{U}(G)$ and thus we have an embedding of
$\mathcal{U}(H)$ into $\mathcal{U}(G)$, and this embedding does not
depend on the choice of $T$.  In fact it will be the same embedding
as described previously.  It follows that
$\mathcal{R}(kH)$ is naturally a subring of $\mathcal{R}(kG)$.
Clearly if $\alpha_1,\dots,
\alpha_n \in \mathcal{U}(H)$ and $t_1, \dots, t_n \in T$, then
$\alpha_1t_1 + \dots + \alpha_nt_n = 0$ if and only if $\alpha_i = 0$
for all $i$, and it follows that if $\beta_1, \dots, \beta_n \in
\mathcal{R}(kH)$ and $\beta_1t_1 + \dots + \beta_nt_n = 0$, then
$\beta_i = 0$ for all $i$.

The above should be compared with the theorem of Hughes
\cite{Hughes70} which we state below.  Recall that a group is
\emph{locally indicable} if every nontrivial
finitely generated subgroup has an infinite cyclic quotient.
Though locally indicable groups are left orderable
\cite[Theorem 7.3.1]{MuraRhemtulla77}
and thus $k*G$ is certainly a domain
whenever $k$ is a division ring, $G$ is a locally indicable group and
$k*G$ is a crossed product, it is still unknown whether such crossed
products can be embedded in a division ring.  Suppose however $k*G$
has a division ring of fractions $D$.  Then we say that $D$ is
\emph{Hughes-free} if whenever $N \lhd H \leqslant G$, $H/N$ is
infinite cyclic, and $h_1, \dots, h_n \in N$ are in distinct cosets
of $N$, then the sum $\mathcal{D}_D(k*N)h_1 + \dots +
\mathcal{D}_D(k*N)h_n$ is direct.
\begin{Thm} \label{THughes}
Let $G$ be a locally indicable group, let $k$ be a division
ring, let $k*G$
be a crossed product, and let $D,E$ be Hughes-free division rings of
fractions for $k*G$.  Then there is an isomorphism
$D \to E$ which is the identity on
$k*G$.
\end{Thm}
This result of Hughes is highly nontrivial, even though the paper
\cite{Hughes70} is only 8 pages long.  This is because the
proof given by Hughes in \cite{Hughes70} is extremely condensed,
and though all the steps are there and correct, it is difficult to
follow.  A much more detailed and somewhat different proof is given
in \cite{DHS03}.

Motivated by Theorem \ref{THughes}, we will extend the definition of
Hughes free to a more general situation.
\begin{Def} \label{DHughes}
Let $D$ be a division ring, let $G$ be a group, let $D*G$ be a
crossed product, and let $Q$ be a ring containing $D*G$ such that
$\mathcal{R}_Q(D*G) = Q$, and every element of $Q$ is either a
zerodivisor or invertible.  In this situation we say that $Q$ is
\emph{strongly Hughes free}
if whenever $N \lhd H \leqslant G$, $h_1, \dots,
h_n \in H$ are in distinct cosets of $N$ and $\alpha_1, \dots,
\alpha_n \in \mathcal{R}_Q(D*N)$, then $\alpha_1h_1 + \dots +
\alpha_nh_n= 0$ implies $\alpha_i= 0$ for all $i$ (i.e.\ the $h_i$
are linearly independent over $\mathcal{R}_Q(D*N)$).
\end{Def}

Then we would like to extend Theorem \ref{THughes} to more general
groups, so we state
\begin{Prob}
Let $D$ be a division ring, let $G$ be a group, let $D*G$ be a
crossed product, and let $Q$ be a ring containing $D*G$ such that
$\mathcal{R}_Q(D*G) = Q$, and every element of $Q$ is either a
zerodivisor or invertible.  Suppose $P,Q$ are strongly Hughes free
rings for $D*G$.  Does there exists an isomorphism $P \to Q$ which is
the identity on $D*G$?
\end{Prob}

It is clear that if $G$ is locally indicable and $Q$ is a division
ring of fractions for $D*G$, then $Q$ is strongly Hughes free implies
$Q$ is Hughes free.  We present the following problem.
\begin{Prob} \label{PHughes}
Let $G$ be a locally indicable group, let $D$ be a division ring, let
$D*G$ be a crossed product, and let $Q$ be a division ring of
fractions for $D*G$ which is Hughes free.  Is $Q$ strongly Hughes
free?
\end{Prob}
It would seem likely that the answer is always ``yes".  Certainly if
$G$ is orderable, then $\mathcal{R}_{D((G))}(D*G)$, the rational
closure (which is the same as the division closure in this case) of
$D*G$ in the Malcev-Neumann power series ring $D((G))$ \cite[Corollary
8.7.6]{Cohn85} is a Hughes free division ring of fractions for
$D*G$.  Therefore by Theorem \ref{THughes} of Hughes, all Hughes free
division ring of fractions for $D*G$ are isomorphic to
$\mathcal{R}_{D((G))}(D*G)$.  It is easy to see that this division
ring of fractions is strongly Hughes free and therefore all Hughes
free division ring of fractions for $D*G$ are strongly Hughes free.

\section{Other Methods}

Embedding $\mathbb{C}G$ into $\mathcal{U}(G)$ has proved to be a very
useful tool, but what about other group rings?  In general we would
like a similar construction when $k$ is a field of nonzero
characteristic.  If $D$ is a division ring, then we can always embed
$D*G$ into a ring in which every element is either a unit or a
zerodivisor, as follows.  Let $V = D*G$ viewed as a right vector
space over $D$, so $V$ has basis $\{\bar{g} \mid g \in G\}$.  Then
$D*G$ acts by left multiplication on $V$ and therefore can be
considered as a subring of the ring of all linear transformations
$\End_D(V)$ of $V$.  This ring is von Neumann
regular.  However it is too large; it is not even directly finite
(that is $xy=1$ implies $yx=1$) when $G$ is infinite.  Another
standard method is to embed $D*G$ in its maximal ring of right
quotients \cite[\S 2.C]{Goodearl76}.  If $R$ is a right
nonsingular ring, then its maximal ring of right quotients
$Q(R)$ is a ring containing $R$ which is a right injective von
Neumann regular ring, and furthermore as a right $R$-module,
$Q(R)$ is the injective hull of $R$
\cite[Corollary 2.31]{Goodearl76}.  By
\cite[Theorem 4]{Snider76}, when $k$ is a field of characteristic
zero, $kG$ is right (and left) nonsingular, consequently $Q(kG)$
is a right self-injective von Neumann regular ring.  However again it
is too large in general.  If $G$ is a nonabelian free group, then
$kG$ is a domain which by Proposition \ref{PLewin}
does not satisfy the Ore condition, so we see from \cite[Exercise
6.B.14]{Goodearl76} that $Q(R)$ is not directly finite.

A very useful technique is that of ultrafilters, see
\cite[p.~76, \S 2.6]{Jacobson80} for example.  We briefly
illustrate this in an example.  Let $k$ be a field and let $G$ be a
group.  Suppose $G$ has a descending chain of normal subgroups
$G=G_0 \geqslant G_1 \geqslant \cdots$ such that $k[G/G_n]$ is
embeddable in a division ring for all $n$.  Then can we embed $kG$ in
a division ring?  It is easy to prove that $kG$ is a domain, but to
prove the stronger statement that $G$ can be embedded in division
ring seems to require the theory of ultrafilters.  For most
applications (or at least for what we are interested in), it is
sufficient to consider ultrafilters on the natural numbers $\mathbb
{N} = \{1,2,\dots\}$.  A filter on $\mathbb {N}$ is a subset
$\omega$ of the power set $\mathcal{P}(\mathbb {N})$ of
$\mathbb{N}$ such that if $X,Y \in \omega$ and
$X \subseteq Z \subseteq \mathbb{N}$, then $X
\cap Y \in \omega$ and $Z \in \omega$.  A filter is proper if it does
not contain the empty set $\emptyset$, and an ultrafilter is a
maximal proper filter.  By considering the maximal ideals in the
Boolean algebra on $\mathcal{P}(\mathbb {N})$, it can be shown that
any proper filter can be embedded in an ultrafilter (this requires
Zorn's lemma), and an ultrafilter has the following properties.
\begin{itemize}
\item
If $X,Y \in \omega$, then $X \cap Y \in \omega$.

\item
If $X \in \omega$ and $X \subseteq Y$, then $Y \in \omega$.

\item
If $X \in \mathcal{P}(\mathbb{N})$, then either $X$
or its complement are in $\omega$.

\item
$\emptyset \notin \omega$.
\end{itemize}
An easy example of an ultrafilter is the set of all subsets
containing $n$ for some
fixed $n \in \mathbb{N}$; such an ultrafilter is called a
\emph{principal ultrafilter}.
An ultrafilter not of this form is called a non-principal
ultrafilter.

Given division rings $D_n$ for $n \in \mathbb {N}$ and
an ultrafilter $\omega$ on $\mathbb {N}$, we can define an equivalence
relation $\sim$ on $\prod_n D_n$ by $(d_1,d_2, \dots) \sim
(e_1,e_2,\dots)$ if and only if there exists $S \in \omega$ such that
$d_n = e_n$ for all $n \in S$.  Then the set of equivalence classes
$(\prod_n D_n)/\sim$ is called the ultraproduct of the division rings
$D_i$ with respect to the ultrafilter $\omega$, and
is a division ring \cite[p.~76, Proposition
2.1]{Jacobson80}.  This can be applied when $R$ is a ring with a
descending sequence of ideals $I_1 \supseteq I_2 \supseteq \dots$ such
that $\bigcap_n I_n = 0$ and $R/I_n$ is a division ring.  The set of
all cofinite subsets of $\mathbb {N}$ is a filter,
so here we let $\omega$ be
any ultrafilter containing this filter.  The corresponding
ultraproduct $D$ of the division rings $R/I_n$ is a division ring.
Furthermore the natural embedding of $R$ into $\prod_n R/I_n$ defined
by $r \mapsto (r+R/I_1,r+R/I_2, \dots)$ induces an embedding of $R$
into $D$.  This proves that $R$ can be embedded in a division ring.

In their very recent preprint \cite{ElekSzabo03},
G\'abor Elek and Endre Szab\'o use
these ideas to embed the group algebra $kG$ over an arbitrary
division ring $k$ in a nice von Neumann regular ring for the class of
\emph{sofic} groups.  The class of sofic groups is a large class of
groups which contains all residually amenable groups and is closed
under taking free products.

Suppose $\{a_n \mid n \in \mathbb {N}\}$ is a bounded sequence of
real numbers and $\omega$ is a non-principal ultrafilter.  Then there
is a unique real number $l$ with the property that given $\epsilon >
0$, then $l$ is in the closure of $\{a_n \mid n \in S\}$ for all $S
\in \omega$.  We call this the
$\omega$-limit of $\{a_n\}$ and write $l = \lim_{\omega} a_n$.

Now let $G$ be a countable amenable group.  Then $G$ satisfies
the F\o lner condition and therefore there exist
finite subsets $X_i$ of $G$ ($i \in \mathbb {N}$) such that
\begin{itemize}
\item
$\bigcup_i X_i = G$.
\item $|X_i| < |X_{i+1}|$ for all $i \in \mathbb {N}$.
\item If $g \in G$, then $\lim_{i\to \infty} |gX_i \cap X_i|/|X_i| =
1$.
\end{itemize}
Let $k$ be a division ring and let $V_i$ denote the right $k$-vector
space with basis $X_i$ ($i \in \mathbb {N}$).  The general element
of $\prod_i \End_k(V_i)$ (Cartesian product) is of the form
$\bigoplus_i \alpha_i$ where $\alpha_i \in \End_k(V_i)$ for all $i$.
For $\beta \in \End_k(V_i)$, we define
\[
\rk_i (\beta) = \frac{\dim_k (\beta V_i)}{\dim_k V_i},
\]
a real number in $[0,1]$.
Now choose a non-principal ultrafilter $\omega$ for $\mathbb {N}$.
Then for $\alpha \in \End_k(V_i)$, we define $\rk(\alpha) =
\lim_{\omega} \rk_n(\alpha_n)$ and $I =
\{\alpha \in \prod_i \End_k(V_i) \mid \rk(\alpha) = 0\}$.  It is
not difficult to check that $I$ is a two-sided ideal of $\prod_i
\End_k(V_i)$.  Now set
\[
\R_k(G) = \frac{\prod_i \End_k(V_i)}{I}
\]
and let $[\alpha]$ denote the image of $\alpha$ in $\R_k(G)$.
Since $\End_k(V_i)$ is von Neumann regular and direct products of von
Neumann regular rings are von Neumann regular, we see that
$\prod_i \End_k(V_i)$ is von Neumann regular and we deduce that
$\R_k(G)$ is also von Neumann regular.
Next we define $\rk([\alpha]) = \rk(\alpha)$.  It can be shown that
$\rk$ is a well-defined rank function \cite[p.~226, Chapter
16]{Goodearl91} and therefore $\R_k(G)$ is directly finite
\cite[Proposition 16.11]{Goodearl91}.

For $g \in G$ and $x \in
X_i$, we can define $\phi(g)x=gx$ if $gx \in X_i$ and
$\phi(g)x = x$ if $gx \notin X_i$.  This determines an embedding
(which is not a homomorphism) of $G$ into
$\prod_i \End_k(V_i)$, and it is shown in \cite{ElekSzabo03}
that the composition with the natural epimorphism $\prod_i
\End_k(V_i)
\twoheadrightarrow \R_k(G)$ yields a homomorphism
$G \to \R_k(G)$.  This homomorphism extends to a ring
homomorphism $\theta \colon kG \to \R_k(G)$ and \cite{ElekSzabo03}
shows that $\ker\theta = 0$.  Thus we have embedded $kG$ into
$\R_k(G)$; in particular this shows that $kG$ is directly
finite because $\R_i(G)$
is.  In fact this construction for $G$ amenable can be extended to
the case $G$
is a sofic group, consequently $kG$ is directly finite if $k$ is a
division ring and $G$ is sofic.  The direct finiteness of $k*G$ for
$k$ a division ring and $G$ free-by-amenable had earlier
been established in \cite{AOP03}.

Another type of localization is considered in \cite{Picavet03}.
Recall that a monoid $M$ is a semigroup with identity, that is $M$
satisfies the axioms for a group except for the existence of
inverses.  If $A$ is a monoid with identity 1, then $M$ is an
$A$-monoid means that there is an action of $A$ on $M$ satisfying
$a(bm) = (ab)m$ and $1m = m$
for all $a,b \in A$ and $m \in M$.  In the case $A$ is a ring
with identity 1 (so $A$ is a monoid under multiplication) and $M$ is
a left $A$-module, then $M$ is an $A$-monoid.  Let $\End(M)$ denote
the monoid of all endomorphisms of the $A$-monoid $M$.  Given a
submonoid $S$ of $\End(M)$, Picavet constructs an
$A$-monoid $S^{-1}M$ with the property that every endomorphism in $S$
becomes an automorphism of $M$, in other words the elements of $S$
become invertible.  To achieve this, he requires that $S$ is a
\emph{localizable} submonoid of End$(M)$.  This means that the
following Ore type conditions hold:

\begin{itemize}
\item
For all $u,v \in S$, there exist $u',v' \in S$ such that $u'u =v'
v$.
\item
For all $u,v,w \in S$ such that $uw=vw$, there is $s \in
S$ such that $su = sv$.
\end{itemize}

The construction is similar to Ore localization.  We describe this
in the case $R$ is ring, $M$ is an $R$-module and $S =
\{\theta^n \mid n \in \mathbb {N}\}$ where $\theta$ is an
endomorphism of $M$.  Clearly $S$ is localizable.  For $m,n \in
\mathbb {N}$ with $m\le n$, we set $M_n = M$ and $\theta_{mn} =
\theta^{n-m} \colon M_m \to M_n$.  Then $(M_n,\theta_{mn})$ forms a
direct system of $R$-modules, and $S^{-1}M$ is the direct limit of
this system.  Clearly $\theta^n$ induces an $R$-automorphism on
$S^{-1}M$ for all $n$, so we have inverted $\theta$.  In the case $R$
is a division ring, $M$ is finitely generated
and $\theta$ is a noninvertible nonnilpotent
endomorphism of $M$, the sequence of $R$-modules $M\theta^n$
eventually stabilizes to a proper nonzero $R$-submodule of $M$, which
is $S^{-1}M$.  It would be interesting to see if this construction
has applications to group rings.

\bibliographystyle{plain}

\end{document}